\documentclass{amsart}
\usepackage{amssymb}
\usepackage{enumerate}

\newtheorem{theorem}{Theorem}
\newtheorem{lemma}[theorem]{Lemma}

\theoremstyle{definition}
\newtheorem{definition}[theorem]{Definition}

\theoremstyle{remark}

\newcommand{\N}{\mathbb{N}}

\DeclareMathOperator{\sper}{Sper}

\DeclareMathOperator{\spec}{Spec}

\newcommand{\p}{\mathfrak{p}}
\newcommand{\q}{\mathfrak{q}}
\newcommand{\du}{\bigcup_{\rightarrow}}
\DeclareMathOperator{\qf}{qf}
\DeclareMathOperator{\Quot}{Quot}

\newcommand{\aff}{f.g.}

\newcommand{\aaff}{f.f.}

\newcommand{\aarch}{almost archimedian}

\DeclareMathOperator{\Nil}{Nil}

\begin{document}
\title[Erratum]{Iterated rings of bounded elements: Erratum}
\author{Markus Schweighofer}
\address{Universit\"at Konstanz\\
         Fachbereich Mathematik und Statistik\\
         78457 Konstanz\\
         Allemagne}
\email{Markus.Schweighofer@uni-konstanz.de}
\thanks{Supported by the project ``Darstellung positiver Polynome''
(Deutsche Forschungsgemeinschaft, Kennung 214371).}
\date{}
\begin{abstract}
We close a gap in the author's thesis \cite{s1,s2}.
\end{abstract}
\maketitle

{\bf The author's proof of \cite[Lemma 4.10]{s1,s2} is not correct.
In this note, we show that this does not affect the validity of any other
statement in \cite{s1,s2}.} We will observe that the lemma in question
holds in any of the following important special cases:
\begin{enumerate}[(a)]
\item $T=\sum A^2$
\item $A$ contains a f.g. subalgebra $C$ such that $T$ is as a preordering
      generated by $T\cap C$.
\item $T$ is as a preordering finitely generated (this is just a special
      case of (b)).
\item $A$ is a reduced ring.
\end{enumerate}
Unfortunately, we don't know whether the lemma holds without any such
additional hypothesis.

\section{Acknowledgements}

The author would like to thank Claus Scheiderer and Mark Olschok who discovered the
gap in a research seminar in Duisburg and reported it to us. The author was able to
close the gap, but it was again Claus Scheiderer who came up with a more
conceptual proof which we will present in Section \ref{redsec} below.

\section{The error}

Recall the situation in the proof of the lemma. We have an extension
$B\subseteq A$ of preordered
rings, i.e., an extension $B\subseteq A$ of rings such that $A$ is equipped
with a preordering $T$ and $B$ with $T\cap B$. We have
$\p\in\spec A$ and $\q=\p\cap B$. Then $B/\q$ can be
viewed as a subring of $A/\p$ but perhaps {\bf not} as a {\bf preordered}
subring, contrary to what is said in \cite{s1,s2}. Here $A/\p$ and $B/\q$
are equipped with the preorderings
$$T_{A/\p}=\{t+\p\mid t\in T\}\qquad\text{and}\qquad
  T_{B/\q}=\{t+\q\mid t\in T\cap B\},$$
respectively. The problem is that for some $b\in B$ it could happen that
$b-t$ lies in $\p$ for some $t\in T$ but not in $\q$ for any $t\in T$.
Consequently, it is not guaranteed whether the preorderings on the quotient
field $\qf(A/\p)=\qf(B/\q)$ induced (or generated) by $T_{A/\p}$ and
$T_{B/\q}$ coincide. Hence it is not clear whether $T\subseteq P$ when
$P$ is chosen like in the proof under review.

\section{Cases where the proof still works}

In case (a), $T\subseteq P$ holds trivially.
In case (b), we may assume that $C\subseteq B$. This implies
$T\cap C\subseteq T\cap B\subseteq Q\subseteq P$ and therefore $T\subseteq P$.

\section{Proof of the lemma for reduced rings}\label{redsec}

In this section, we will prove that \cite[Lemma 4.10]{s1,s2}
holds under the additional assumption (d) that $A$ is a \emph{reduced}
ring, i.e., contains no nonzero nilpotent elements.
This is done in Lemma \ref{nuages} below.

We need some well-known facts from commutative ring theory whose proofs we
include for the convenience of the reader. Let $A$ always denote a
commutative ring (with unity, of course).

\begin{lemma}\label{triv1}
Let $\p_1,\dots,\p_n$ be prime ideals of a commutative ring $A$ satisfying
$\p_i\not\subseteq\p_j$ for $i\neq j$. Then
$\p_i\not\subseteq\bigcup_{j\neq i}\p_j$ for all $i\in\{1,\dots,n\}.$
\end{lemma}

\begin{proof}
For every $(i,j)$ with $j\neq i$, choose $a_{ij}\in\p_i\setminus\p_j$.
Then we have for instance
$$\sum_{i=2}^n\prod_{j\neq i}a_{ij}\in\p_1\setminus
  \bigcup_{j\neq i}\p_j.$$
\end{proof}

\begin{lemma}\label{triv2}
Let $A$ be a reduced ring with only finitely many pairwise distinct
minimal prime ideals
$\p_1,\dots,\p_n$. Then the zero divisors in $A$ (i.e., the elements $a\in A$
for which there is some $0\neq b\in A$ with $ab=0$) are exactly the elements
of $\p_1\cup\ldots\cup\p_n$.
\end{lemma}

\begin{proof}
Consider an element of $\p_1\cup\ldots\cup\p_n$, say $a\in\p_1$.
By the preceding lemma, we can choose $b_i\in\p_i\setminus\bigcup_{j\neq i}
\p_j$ for $i\in\{2,\dots,n\}$. Then $ab_2\cdots b_n\in\p_1\cap\ldots\cap\p_n=0$
but $b_2\cdots b_n\notin\p_1$, in particular $b_2\cdots b_n\neq 0$.
Thus $a$ is a zero divisor.

Conversely, suppose $a\in A\setminus (\p_1\cup\ldots\cup\p_n)$. We show that
$a$ is not a zero divisor. Suppose therefore that $b\in A$ and $ab=0$.
From $ab=0\in\p_i$ and $a\notin\p_i$ it follows that $b\in\p_i$ for all
$i\in\{1,\dots,n\}$. Hence $b$ lies in $\p_1\cap\ldots\cap\p_n=0$.
\end{proof}

\begin{definition}\label{triv3}
The \emph{total quotient ring} of $A$ is the ring
$$\Quot(A):=S^{-1}A$$
where $S$ is the set non zero divisors of $A$ (note that $1\in S$ and
$SS\subseteq S$).
\end{definition}

\begin{lemma}\label{triv4}
The canonical homomorphism $A\to\Quot(A)$ is an embedding.
\end{lemma}

\begin{proof}
Suppose $a\in A$ and $a/1=0\in\Quot(A)$. Then there is some non zero
divisor $s$ of $A$ such that $as=0$ in $A$. But then, of course, $a=0$.
\end{proof}

\begin{lemma}\label{triv5}
Let $A$ be reduced with only finitely many pairwise distinct minimal prime
ideals $\p_1,\dots,\p_n$. Then there is a canonical isomorphism
$$\Quot(A)\xrightarrow\cong\qf(A/\p_1)\times\dots\times\qf(A/\p_n).$$
\end{lemma}

\begin{proof}
Let $S$ denote the set of non zero divisors of $A$. Then $\Quot(A)=S^{-1}A$
by Definition \ref{triv3}. The prime ideals of $S^{-1}A$ correspond to the
prime ideals of $A$ contained in $A\setminus S=\p_1\cup\ldots\cup\p_n$
(Lemma \ref{triv2}). But by Lemma \ref{triv1}, the only prime
ideals of $A$ contained in $\p_1\cup\ldots\cup\p_n$ are $\p_1,\dots,\p_n$
themselves. Therefore $S^{-1}\p_1,\dots,S^{-1}\p_n$ are (all) pairwise distinct
maximal ideals of $S^{-1}A$. In particular, these ideals are pairwise coprime.
Moreover it is easy to see that $$S^{-1}\p_1\cap\ldots\cap S^{-1}\p_n=0$$
and that there is a canonical isomorphism
$$S^{-1}A/S^{-1}\p_i\xrightarrow\cong\qf(A/\p_i)$$
for all $i\in\{1,\dots,n\}$. Our claim therefore follows
by applying the Chinese Remainder Theorem to the ring $S^{-1}A$.
\end{proof}

\begin{lemma}\label{triv6}
Let $A$ be reduced with only finitely many pairwise distinct minimal prime
ideals $\p_1,\dots,\p_n$. Let $B$ be a subring of $A$ such that the following
conditions hold:
\begin{enumerate}[(1)]
\item\label{localcond}
The canonical embedding $\qf(B/(\p_i\cap B))\hookrightarrow\qf(A/\p_i)$
is a field isomorphism for all $i$.
\item\label{mutcond}
$\p_i\cap B\not\subseteq\p_j\cap B$ for all $i\neq j$.
\end{enumerate}
Then there is a canonical isomorphism $\Quot(B)\xrightarrow\cong\Quot(A)$.
\end{lemma}

\begin{proof}
Every minimal prime ideal of $B$ is of the form $\p_i\cap B$ for some
$i\in\{1,\dots,n\}$ (confer \cite[Remark 4.8]{s1,s2}). Conversely, we argue
that every $\p_i\cap B$ is actually a minimal prime ideal of $B$. To see this,
observe that $\p_i\cap B$ contains in any case a minimal prime ideal. Hence
$\p_j\cap B\subseteq\p_i\cap B$ for some minimal prime ideal $\p_j\cap B$
of $B$. Condition (\ref{mutcond}) forces $i=j$ showing that $\p_i\cap B$ is
itself a minimal prime ideal.
Our claim follows now from condition (\ref{localcond}), the preceding
lemma and the fact that $\p_1\cap B,\dots,\p_n\cap B$ are exactly the
\emph{pairwise distinct} minimal prime ideals of $B$.
\end{proof}

\begin{lemma} \label{nuages}
Suppose $A$ is \aaff , \aarch\ and reduced. Then
\[A=\du B,\quad\text{where $B$ ranges over \aff , \aarch\ algebras.}\]
\end{lemma}

\begin{proof}
Denote the finitely many minimal pairwise distinct prime ideals of $A$ by
$\p_1,\dots,\p_n$. Since $A$ is f.f., we can choose $a_1,\dots,a_m\in A$ such
that $\qf(A/\p_i)$ is for each $i\in\{1,\dots,n\}$ generated by
$a_1+\p_i,\dots,a_m+\p_i$ as a field over $K$. Choose moreover
$b_{ij}\in\p_i\setminus\p_j$ for all $i\neq j$. Clearly, $A=\du B$ where $B$
ranges over all f.g. subalgebras of $A$ containing all of the finitely many
$a_i$ and $b_{ij}$.

Now we fix such a $B$. It remains to show that $B$ is almost archimedean.
Fix an arbitrary $Q\in\sper B$ such that $Q\cap -Q$ is a minimal prime ideal
of $B$. We
have to show that $Q$ is archimedean. There is $i$ such that
$Q\cap -Q=\p_i\cap B$ by \cite[Remark 4.8]{s1,s2}. Since $B$ contains all
$a_i$ and
$b_{ij}$, conditions (\ref{localcond}) and (\ref{mutcond}) from Lemma
\ref{triv6} are satisfied. From (\ref{localcond}) we see that
$\qf(A/\p_i)=\qf(B/(\p_i\cap B))$ as fields (not necessarily as preordered
fields!).
Consequently, there is some ordering $P$ of the ring $A$ such that
$Q=P\cap B$ and $P\cap -P=\p_i$. It remains to show that
$P\in\sper A$ or, in more explicit words, $T\subseteq P$. Once we have shown
this, it follows that $P$ is archimedean since $A$ is almost archimedean.
But then $Q$ must of course be archimedean, too.

To show $T\subseteq P$, we use Lemma \ref{triv5} saying that there is a
canonical isomorphism $\Quot(B)\xrightarrow\cong\Quot(A)$.
Consider an arbitrary $t\in T$. Since $t/1$ lies in the image of this
isomorphism,
there is some $b\in B$ and some non zero divisor $s$ in $B$ such that
$b/s=t/1$ holds in $\Quot(A)$. This implies $b/1=st/1$ in $\Quot(A)$ and a
fortiori $b=st$ in $A$ by Lemma \ref{triv4}. Hence
$s^2t\in T\cap B\subseteq Q\subseteq P$.
If $s$ were an element of $P\cap -P=\p_i$, then it would lie in the
minimal prime ideal $\p_i\cap B=Q\cap -Q$ of $B$ which is impossible by Lemma
\ref{triv2}. From $s\not\in P\cap -P$ it follows now that $t\in T$.
\end{proof}

\section{Closing the gap}

Finally, we show that Lemma \ref{nuages} is enough to ensure the validity
of all results of \cite{s1,s2}, with the only possible exception of
\cite[Lemma 4.10]{s1,s2}, of course. In fact, \cite[Lemma 4.10]{s1,s2} is
only applied once in \cite{s1,s2}, namely in the proof of
\cite[Theorem 4.13]{s1,s2} (and its variant for quadratic modules
and semiorderings in \cite[Subsection 6.2]{s1,s2} which can be treated
completely analogously).

Hence it suffices to show that we can restrict ourselves in the proof of
\cite[Theorem 4.13]{s1,s2} to \emph{reduced} $A$ since then we can apply
Lemma \ref{nuages} instead of \cite[Lemma 4.10]{s1,s2}. Assume therefore
that \cite[Theorem 4.13]{s1,s2} has already been shown for reduced $A$.

Now for general $A$, denote the nilradical of $A$ by $\Nil(A)$. Suppose
$A=H(A)$.
Then $A/\Nil(A)=H(A/\Nil(A))$ and therefore $A/\Nil(A)=H'(A/\Nil(A))$ since
$A/\Nil(A)$ is reduced. Suppose $a\in A$. We show that $a\in H'(A)$.
From $a+\Nil(A)\in A/\Nil(A)=H'(A/\Nil(A))$ we obtain $\nu\in\N$ and
$b\in\Nil(A)$ such that $\nu-a+b\in T$. Clearly $b\in H'(A)$ by
\cite[Lemma 4.1]{s1,s2}. This supplies us with a $\nu'\in\N$ such that
$\nu'-b\in T$. Finally, we see that
$$(\nu+\nu')-a=(\nu-a+b)+(\nu'-b)\in T+T\subseteq T.$$
Since $a\in A$ was arbitrary, we see that $A=H'(A)$ as desired.

\end{document}